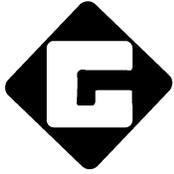

# The roots of any polynomial equation

**G.A.Uytdewilligen**,
*Bergen op Zoomstraat 76, 5652 KE Eindhoven. g.a.uytdewilligen@zonnet.nl*

**Abstract**

We provide a method for solving the roots of the general polynomial equation

$$a_n \cdot x^n + a_{n-1} \cdot x^{n-1} + \ldots + a_1 \cdot x + s = 0 \tag{1}$$

To do so, we express x as a powerseries of s, and calculate the first n-1 coefficients. We turn the polynomial equation into a differential equation that has the roots as solutions. Then we express the powerseries' coefficients in the first n-1 coefficients. Then the variable s is set to a0. A free parameter is added to make the series convergent. © 2004 G.A.Uytdewilligen. All rights reserved.

**The method**

The method is based on [1]. Let's take the first n-1 derivatives of (1) to s. Equate these derivatives to zero. Then find $\frac{d^i}{ds^i} x(s)$ in terms of x(s) for i from 1 to n-1. Now make a new differential equation

$$m_1 \cdot \frac{d^{n-1}}{ds^{n-1}} x(s) + m_2 \cdot \frac{d^{n-2}}{ds^{n-2}} x(s) + \ldots + m_n \cdot x(s) + m_{n+1} = 0 \tag{2}$$

and fill in our $\frac{d^i}{ds^i} x(s)$ in (2). Multiply by the denominator of the expression. Now we have a polynomial in x(s) of degree higher then n. Using (1) as property, we simplify this polynomial to the degree of n. Set it equal to (1) and solve $m_1 \ldots m_{n+1}$ in terms of s and $a_1 \ldots a_n$ Substituting these in (2) gives a differential equation that has the zeros of (1) among its solutions. We then insert

$$x(s) = y(s) - \frac{a_{n-1}}{n \cdot a_n} \tag{3}$$

in (2). Multiplying by the denominator we get a differential equation of the linear form:

$$p_1 \cdot \frac{d^{n-1}}{ds^{n-1}} y(s) + p_2 \cdot \frac{d^{n-2}}{ds^{n-2}} y(s) + \ldots + p_n \cdot y(s) = 0 \tag{4}$$

With $p_1(s)..p_n(s)$ polynomials in s. If we substitute our powerseries, all the coefficients are determined by the first n-1 coefficients. The first coefficients are calculated as follows: A powerseries is filled in in (1).

$$x(s) = \sum_{i=0}^{n-2} b_i \cdot s^i - \frac{a_n}{n \cdot a_{n-1}} \tag{5}$$

and it should be zero for all s. From this, we calculate $b_i$ for i from 0 to n-2. $b_0$ Is a root af an n-1 degree polynomial and the other $b_i$ are expressed in $b_0$ Now a powerseries is inserted in (4):

$$y(s) = \sum_{i=0}^{\infty} b_i \cdot s^i$$

(6)

and we get an equation of the form:

$$q_1(i) \cdot c_1 \cdot b_i + q_2(i) \cdot c_2 \cdot b_{i+1} + \ldots + q_n(i) \cdot c_n \cdot b_{i+n-1} = 0$$

(7)

where $q_m(i)$ are polynomials in i of degree n-1. $c_m$ Are constants.
We define $b_{n-1}$ as the determinant of a matrix A

$$A = \begin{vmatrix} \frac{c_{n-1} \cdot q_{n-1}(0)}{c_n \cdot q_n(0)} & .. & \frac{c_2 \cdot q_2(0)}{c_n \cdot q_n(0)} & \frac{c_1 \cdot q_1(0)}{c_n \cdot q_n(0)} & 0 & 0 & 0 & .. & 0 & 0 \\ 0 & .. & 0 & 0 & 1 & 1 & 0 & .. & 0 & 0 \\ 0 & .. & 0 & 0 & 0 & 1 & 1 & .. & 0 & 0 \\ .. & .. & .. & .. & .. & .. & .. & .. & .. & .. \\ 0 & .. & 0 & 0 & 0 & 0 & 0 & .. & 1 & 1 \\ 0 & .. & 0 & 1 & b_0 & 0 & 0 & .. & 0 & 0 \\ 0 & .. & 1 & 0 & 0 & -b_1 & 0 & .. & 0 & 0 \\ 0 & .. & 0 & 0 & 0 & 0 & b2 & .. & 0 & 0 \\ .. & .. & .. & .. & .. & .. & .. & .. & .. & .. \\ 1 & 0 & 0 & 0 & 0 & 0 & 0 & .. & 0 & -1^{n-2} \cdot b_{n-2} \end{vmatrix}$$

(8)

and for the rest of the coefficients

$$b_{i+n-1} = \begin{vmatrix} \frac{-c_{n-1} \cdot q_{n-1}(i)}{c_n \cdot q_n(i)} & \frac{-c_{n-2} \cdot q_{n-2}(i)}{c_n \cdot q_n(i)} & \frac{-c_{n-3} \cdot q_{n-3}(i)}{c_n \cdot q_n(i)} & . & \frac{-c_1 \cdot q_1(i)}{c_n \cdot q_n(i)} & 0 & 0 & 0 & 0 & .. \\ -1 & \frac{-c_{n-1} \cdot q_{n-1}(i-1)}{c_n \cdot q_n(i-1)} & \frac{-c_{n-2} \cdot q_{n-2}(i-1)}{c_n \cdot q_n(i-1)} & . & \frac{-c_2 \cdot q_2(i-1)}{c_n \cdot q_n(i-1)} & \frac{c_1 \cdot q_1(i-1)}{c_n \cdot q_n(i-1)} & 0 & 0 & 0 & .. \\ .. & .. & .. & .. & .. & .. & .. & .. & .. & .. \\ 0 & 0 & 0 & -1 & \frac{-c_{n-1} \cdot q_{n-1}(1)}{c_n \cdot q_n(1)} & \frac{c_{n-2} \cdot q_{n-2}(1)}{c_n \cdot q_n(1)} & . & \frac{c_1 \cdot q_1(1)}{c_n \cdot q_n(1)} & 0 & .. \\ 0 & 0 & 0 & 0 & -1 & & & & & \\ 0 & 0 & 0 & 0 & 0 & & & & & \\ 0 & 0 & 0 & 0 & 0 & & & A & & \\ 0 & 0 & 0 & 0 & 0 & & & & & \\ 0 & 0 & 0 & 0 & 0 & & & & & \\ .. & .. & .. & .. & .. & & & & & \end{vmatrix}$$

(9)

The series (6) can be proven to be convergent [2] if for a constant E

$$s \leq \frac{1}{E} \cdot |c_n|$$

(10)

and if all the absolute values of the coefficients of (1) are chosen smaller then 1. This is done by dividing a polynomial by (more than) the maximum of the absolute values of the coefficients.

To make the series convergent, we transform s to e·s. That is, if we insert the powerseries it is not in s but in e·s. Writing out the terms of the sum, we find that each term $d_i$ has a factor $s^i \cdot e^{i-n+1}$. Setting

$$e = \frac{|c_n|}{E} \qquad (12)$$

We still need s<1, which is why we set s to $a_0$ and $a_0$<1.